\documentclass[12pt]{amsart}

\usepackage{amsmath,amssymb,amsthm,amscd}
\usepackage{verbatim}
\usepackage{comment}
\usepackage{multirow}
\usepackage{mathtools}
\usepackage[shortlabels]{enumitem}
\usepackage[margin=1in]{geometry} 
\usepackage{graphics}
\usepackage{pgfplots}
\pgfplotsset{compat=1.15}
\usepackage{mathrsfs}
\usetikzlibrary{arrows}
\numberwithin{equation}{section}
\usepackage{caption} 
\usepackage{longtable}
\usepackage{ltablex}
\usepackage{float}

\newtheorem{theorem}{Theorem}[section]

\newtheorem{lemma}[theorem]{Lemma}

\newtheorem{corollary}[theorem]{Corollary}

\newtheorem{question}[theorem]{Question}

\newtheorem*{theorem*}{Theorem}

\theoremstyle{definition}
\newtheorem{definition}[theorem]{Definition}

\newtheorem{warning}[theorem]{Warning}

\newtheorem{remark}[theorem]{Remark}

\newcommand{\PP}{ \ensuremath{\mathbb{P}}}

\def\cocoa{{\hbox{\rm C\kern-.13em o\kern-.07em C\kern-.13em o\kern-.15em A}}}

\newcommand{\specialcell}[2][c]{%
\renewcommand*{\arraystretch}{1}
  \begin{tabular}[#1]{@{}c@{}}#2\end{tabular}}

\begin{document}

\title{Geproci sets on skew lines in $\mathbb P^3$ with two transversals}

\author[L.~Chiantini et al.]{Luca Chiantini}
\address[L.~Chiantini]{Dipartmento di Ingegneria dell'Informazione e Scienze Matematiche, Universit\`a di Siena, Italy}
\email{luca.chiantini@unisi.it}

\author[]{Pietro De Poi}
\address[P.~De Poi]{Dipartimento di Scienze Matematiche, Informatiche e Fisiche, Universit\`a degli Studi di Udine, 
Via delle Scienze, 206, 33100 Udine, Italy}
\email{pietro.depoi@uniud.it}

\author[]{{\L}ucja Farnik}
\address[{\L}.~Farnik]{Department of Mathematics, University of the National Education Commission,
   Podcho\-r\c a\.zych~2,
   PL-30-084 Krak\'ow, Poland}
\email{lucja.farnik@gmail.com}

\author[]{Giuseppe Favacchio}
\address[G.~Favacchio]{Dipartimento di Ingegneria, Universit\`a degli studi di Palermo,
Viale delle Scienze,  90128 Palermo, Italy}
\email{giuseppe.favacchio@unipa.it}

\author[]{Brian Harbourne}
\address[B.~Harbourne]{Department of Mathematics,
University of Nebraska,
Lincoln, NE 68588-0130 USA}
\email{brianharbourne@unl.edu}

\author[]{Giovanna Ilardi}
\address[G.~Ilardi]{Dipartimento di Matematica e Applicazioni ``Renato Caccioppoli'', 
Universit\`a degli Studi Napoli ``Federico II'', Via Cintia, Monte S. Angelo, 80126 Napoli, Italy}
\email{giovanna.ilardi@unina.it}

\author[]{Juan Migliore} 
\address[J.~Migliore]{Department of Mathematics,
University of Notre Dame,
Notre Dame, IN 46556 USA}
\email{migliore.1@nd.edu}

\author[]{Tomasz Szemberg}
\address[T.~Szemberg]{Department of Mathematics, University of the National Education Commission,
   Podcho\-r\c a\.zych~2,
   PL-30-084 Krak\'ow, Poland}
\email{tomasz.szemberg@gmail.com}

\author[]{Justyna Szpond}
\address[J.~Szpond]{Department of Mathematics, University of the National Education Commission,
   Podcho\-r\c a\.zych~2,
   PL-30-084 Krak\'ow, Poland}
\email{szpond@gmail.com}

\thanks{Chiantini, De Poi, Favacchio and Ilardi are members of the Italian GNSAGA-INDAM}
\thanks{Farnik was partially supported by National Science Centre, Poland, grant 2018/28/C/ST1/00339.
}
\thanks{Favacchio was partially supported by Fondo di
Finanziamento
per la Ricerca di
Ateneo, Università degli studi di Palermo}
\thanks{Harbourne was partially supported by Simons Foundation grant \#524858.}
\thanks{Migliore was partially supported by Simons Foundation grant \#839618.
}
\thanks{Szemberg and Szpond were partially supported by National Science Centre, Poland, grant 2019/35/B/ST1/00723.}


\maketitle

%\tableofcontents

%%%%%%%%%%%%%%%%%%%%%%%%%%%%%%%%%%%%%%%%%%%%%%%%%%%%%%%%%%%
\begin{abstract}
    The purpose of this work is to pursue classification of geproci sets. Specifically we classify $[m,n]$-geproci sets which consist of $m=4$ points on each of $n$ skew lines, assuming the skew lines have two transversals in common. We show that in this case $n\leq 6$. Moreover we show that all geproci sets of this type are contained in the \emph{standard construction} for $m=4$ introduced in  \cite{POLITUS}. Finally, we propose a conjectural representation for all geproci sets of this type, irrespective of the number $m$ of points on each skew line. 
\end{abstract}

%%%%%%%%%%%%%%%%%%%%%%%%%%%%%%%%%%%%%%%%%%%%%%%%%%%%%%%%%%%
\section{Introduction}

Throughout this paper we work over the complex numbers, and $Z$ will always be a reduced finite set of points in $\mathbb P^3$.  
We denote by $\overline{Z}_{P,H}$ (but often just by $\overline{Z}$) 
the image of $Z$ under projection to a  plane $H\cong\mathbb P^2$ 
from a general point $P$. 
When $\overline{Z}_{P,H}$ is a transverse intersection of two curve in $H$ we say $Z$ is {\it geproci}.

If $Z$ is itself a complete transverse intersection of two curves in a plane, then it is easy to see that $\overline{Z}_{P,H}$ will  be
a complete transverse intersection of two curves in $H$, and hence that $Z$ is geproci. The question of whether 
nondegenerate (i.e., non-coplanar) examples of geproci sets $Z$ exist was raised by Polizzi and answered by Panov \cite{Polizzi, CDFPGR} who pointed to
grids (i.e., intersections $Z$ of two curves $A$ and $B$, each consisting of skew lines, such that every component of $A$ meets every component of $B$ transversely). 

Whether nondegenerate non-grid geproci examples existed remained open until it was noticed (see \cite[Appendix]{CM}), based on recent work on
unexpectedness \cite{HMNT}, that certain root systems (such as $D_4$ and $F_4$) gave examples 
of nondegenerate non-grid geproci sets.
These examples were key in identifying
a broad class of examples called half grids
\cite{PSS} which are our main focus here.

More specifically, following \cite{POLITUS}, we say that $Z$ is an {\em $(a,b)$-geproci set} if $\overline{Z}$
is the transverse intersection of curves in $H$ of degrees $a$ and $b$ with $a \leq b$;
i.e., if $\overline{Z}$ is a complete intersection of type $(a,b)$ with $a \leq b$,
and we say that $Z$ is {\em  $\{a,b\}$-geproci} if we drop the condition $a \leq b$.

\begin{definition}\label{def:grid and half grid}
An $(a,b)$-geproci set is an $(a,b)$-{\it grid} if 
there is a set $A$ of
$a$ skew lines with each line containing exactly $b$ of the points, and a set $B$ of
$b$ skew lines with each line containing exactly $a$ of the points (if $a=b$ we also require $A\cap B=\varnothing$; this is automatic if $a<b$).
An $(a,b)$-{\it half grid} 
(or $\{a,b\}$-{\it half grid})
is an $(a,b)$-geproci (or $\{a,b\}$-geproci, resp.) set for which 
either $A$ or $B$ exists but not both.
% there is either a set of
% $a$ skew lines with each line containing $b$ % of the points, or a set of
% $b$ skew lines with each line containing $a$ % of the points, but not both.
In addition, we say that an $\{a,b\}$-geproci set is an $[a,b]$-\emph{half grid} if it consists of
$a$ points on each of $b$ skew lines. 
\end{definition}

The main results of \cite{POLITUS} establish the existence of non-grid $(a,b)$-geproci sets of points for all integers $4\leq a\leq b$ and for $(a,b)=(3,4)$. In the latter case \cite{POLITUS} provides also the full classification: the only non-grid $(3,4)$-geproci set in $\PP^3$ is determined by the $D_4$ root system. This result has a profound impact on the present note.

The next natural case to study are half grid $(4,4)$-geproci sets. They were fully classified in \cite{POLITUS2}, where the authors show that there are only two possible cases. 
%{\color{blue} Is this true? I thought \cite{POLITUS2} {\em assumed} half grid. Where is it shown that geproci implies half grid for $(4,4)$?}

Moreover all but three geproci sets found up to now are half grids; see Definition \ref{def:grid and half grid} above. Working under this assumption, we extend the detailed classification of geproci sets to $[4,n]$-half grids for $n\geq 4$. Our main result is the following.
\begin{theorem}\label{thm:main}
Let $Z$ be a $(4,n)$-half grid of 4 points on each of $n$ lines such that there are two lines transversal to the $n$ half grid lines. Then $n\leq 6$ and $Z$ is projectively equivalent to a subset of the $F_4$ configuration.
\end{theorem}

\section{Preliminaries}
Here we recall some basic notions and facts we shall use in the sequel. We begin with the following useful observation, which  is a direct consequence of the classification of all $(3,4)$-geproci sets performed in \cite{POLITUS}.

\begin{lemma}\label{3lines}
    Let $Z$ be a $[4,n]$-half grid with $n\geq 4$. Then any subset $W\subset Z$ consisting of all points in $Z$ on $3$ of the half grid lines is a $(3,n)$-grid.
\end{lemma}

Next we recall two basic notions from projective geometry.

\begin{definition} 
Recall that the \textit{cross ratio} of an ordered set of four distinct points $P_1=[x_1:y_1],P_2=[x_2:y_2],P_3=[x_3:y_3],P_4=[x_4:y_4]$ with respect to some (in fact: any)  choice of coordinates on $\PP^1$ is
\begin{equation*} 
j(P_1,P_2;P_3,P_4)= \frac {(x_1y_3-y_1x_3)(x_2y_4-y_2x_4)}{(x_1y_4-y_1x_4)(x_2y_3-y_2x_3)}.
\end{equation*}
\end{definition}

\begin{definition}	 We say that  the points are {\em harmonic} if their cross ratio is $-1$, $1/2$ or $2$ (the specific value
depends on the ordering of the points).

We say that the points are {\em anharmonic} if their cross ratio is $1/2 + \sqrt{3}i/2$ or
$1/2 - \sqrt{3}i/2$.
\end{definition}
\begin{warning}
In this note $(a,b,c,d)$ with $\left\{a,b,c,d\right\}=\left\{1,2,3,4\right\}$ denotes a permutation which sends $1$ to $a$, $2$ to $b$, $3$ to $c$ and $4$ to $d$. So this is not the cycle notation!
\end{warning}
\begin{remark}\label{rem:CR invariant}
It is well known and easy to check by direct calculation that for any $4$-tuple of mutually distinct points
$$j(P_1,P_2;P_3,P_4)=
j(P_{\sigma(1)},P_{\sigma(2)};P_{\sigma(3)},P_{\sigma(4)})$$
for $\sigma\in\left\{id, (2,1,4,3), (3,4,1,2), (4,3,2,1)\right\}$. Note that the non-trivial permutations leaving invariant the cross-ratio of an arbitrary set of $4$ points  are exactly fixed point free involutions in $S_4$.\\
For harmonic points $P_1,\ldots,P_4$ the set of permutations leaving them invariant is bigger. In addition to the four permutations mentioned above, it contains also the following four elements:
$$ (2,1,3,4), (1,2,4,3), (3,4,2,1), (4,3,1,2).$$
The first two of these elements are involutions with a fixed point and the other two 
%$(3,4,2,1), (4,3,1,2)$ 
are $4$-cycles.
\end{remark}
We conclude this section with the following well-known and useful observation exploring the cross-ratio; for a proof see, e.g., \cite[Paragraph 3.4.1]{EisHar16}.
\begin{lemma}\label{lem:4 points and automorphis}
Let $P_1,\ldots,P_4$ and $R_1,\ldots,R_4$ be two four-tuples of points on the projective line $\PP^1$. If
$$j(P_1,P_2;P_3,P_4)=j(R_1,R_2;R_3,R_4),$$
then there exists a linear projective map $F:\PP^1\to\PP^1$ such that $F(P_i)=R_i$ for $i=1,\ldots,4$.\end{lemma}

\section{Permutations on half grids}\label{sec: permutations on half grids}

\subsection{Classification of $[4,n]$-half grids with transversals and containing a $(4,4)$-grid}

Let $Z$ be a $[4,n]$-half grid, with two transversals $T_1,T_2$. This means that there are $n$ skew lines $L_1,\ldots,L_n$, each containing exactly $4$  points from $Z$ such that all these lines intersect lines $T_1$ and  $T_2$. 

We assume additionally that none of the intersection points between the lines $L_1,\ldots,L_n$ and $T_1, T_2$ belongs to $Z$. Furthermore, we assume that $Z$ contains a $(4,4)$-grid spanning a smooth quadric $Q$ (note that we do not need to assume that $Z$ does not contain a $(4,5)$-grid -- see Corollary \ref{max4}). 
More specifically, let us suppose that $L_1,\dots,L_4$  are in the ``vertical'' ruling of $Q$ and let us denote by $M_1,\dots,M_4$ the grid lines in the ``horizontal'' ruling.
The points $P_{ij}=L_i\cap M_j$ for $1\leq i,j\leq 4$ form a grid. The transversals $T_1,T_2$ are contained in $Q$ and they are lines in the ``horizontal'' ruling different from $M_1,\ldots,M_4$. 

Since $Z$ is not a grid, it must be $n\geq 5$. For every line $L\in\left\{L_5,\ldots,L_n\right\}$ and $i\in\left\{2,3,4\right\}$ we denote by 
$Q_i^L$ the quadric spanned by $L_1,L_i$ and $L$.
This data determines an element of the symmetric group $S_4$, that we denote with $\sigma^L_i$, as follows. For a point $P_{ij}$ on $L_i$, by Lemma \ref{3lines}, there is a line in $Q^L_i$ in the ruling of the transversals passing through $P_{ij}$ and meeting $L_1$ in a point of $Z$, say $P_{1k}$. We define the permutation $\sigma_i^L$ by putting $\sigma_i^L(j)=k$.

\definecolor{uuuuuu}{rgb}{0,0,0}
\definecolor{xdxdff}{rgb}{0,0,0}
\definecolor{ududff}{rgb}{0,0,0}
\begin{figure}[ht]
    \centering
\begin{tikzpicture}[line cap=round,line join=round,>=triangle 45,x=1cm,y=1cm]
\clip(-1,-0.5) rectangle (8,9);
\draw [line width=2pt] (1,0) -- (1,4.5);
\draw [line width=2pt] (2,0) -- (2,4.5);
\draw [line width=2pt] (3,0) -- (3,4.5);
\draw [line width=2pt] (4,0) -- (4,4.5);
\draw [line width=2pt,domain=0.5:5] plot(\x,{(--1-0*\x)/1});
\draw [line width=2pt,domain=0.5:5] plot(\x,{(--2-0*\x)/1});
\draw [line width=2pt,domain=0.5:5] plot(\x,{(--3-0*\x)/1});
\draw [line width=2pt,domain=0.5:5] plot(\x,{(--4-0*\x)/1});
\draw [line width=2pt,domain=2:7] plot(\x,{(--34.1076-2.44*\x)/3.86});
\draw [line width=2pt,domain=0:5] plot(\x,{(--2--1*\x)/1});
\begin{scriptsize}
\draw [fill=ududff] (1,1) circle (2.5pt);
\draw[color=ududff] (0.5,0.5) node {$P_{1,1}$};
\draw [fill=ududff] (1,2) circle (2.5pt);
\draw[color=ududff] (0.5,1.5) node {$P_{1,2}$};
\draw[color=black] (1,5) node {$L_1$};
\draw [fill=ududff] (2,1) circle (2.5pt);
\draw[color=ududff] (1.6,0.5) node {$P_{2,1}$};
\draw [fill=ududff] (2,2) circle (2.5pt);
\draw[color=ududff] (1.6,1.5) node {$P_{2,2}$};
\draw[color=black] (2,5) node {$L_2$};
\draw [fill=ududff] (3,1) circle (2.5pt);
\draw[color=ududff] (3.3,1.3) node {$P_{3,1}$};
\draw [fill=ududff] (3,2) circle (2.5pt);
\draw[color=ududff] (3.3,2.3) node {$P_{3,2}$};
\draw[color=black] (3,-0.4) node {$L_3$};
\draw [fill=ududff] (4,1) circle (2.5pt);
\draw[color=ududff] (4.3,1.3) node {$P_{4,1}$};
\draw [fill=ududff] (4,2) circle (2.5pt);
\draw[color=ududff] (4.3,2.3) node {$P_{4,2}$};
\draw[color=black] (4,-0.4) node {$L_4$};
\draw[color=black] (5.4,1) node {$M_1$};
\draw[color=black] (5.4,2) node {$M_2$};
\draw [fill=xdxdff] (1,3) circle (2.5pt);
\draw[color=xdxdff] (0.5,3.4) node {$P_{1,3}$};
\draw [fill=xdxdff] (2,3) circle (2.5pt);
\draw[color=xdxdff] (1.6,2.5) node {$P_{2,3}$};
\draw[color=black] (5.4,3) node {$M_3$};
\draw [fill=xdxdff] (1,4) circle (2.5pt);
\draw[color=xdxdff] (0.5,4.4) node {$P_{1,4}$};
\draw [fill=xdxdff] (2,4) circle (2.5pt);
\draw[color=xdxdff] (1.6,4.4) node {$P_{2,4}$};
\draw[color=black] (5.4,4) node {$M_4$};
\draw [fill=uuuuuu] (3,3) circle (2pt);
\draw[color=uuuuuu] (3.3,3.3) node {$P_{3,3}$};
\draw [fill=uuuuuu] (3,4) circle (2pt);
\draw[color=uuuuuu] (3.3,4.3) node {$P_{3,4}$};
\draw [fill=uuuuuu] (4,4) circle (2pt);
\draw[color=uuuuuu] (4.3,4.3) node {$P_{4,4}$};
\draw [fill=uuuuuu] (4,3) circle (2pt);
\draw[color=uuuuuu] (4.3,3.3) node {$P_{4,3}$};
\draw [fill=ududff] (2.62,7.18) circle (2.5pt);
%\draw[color=ududff] (2.78,7.61) node {$Q$};
\draw [fill=ududff] (6.48,4.74) circle (2.5pt);
%\draw[color=ududff] (6.64,5.17) node {$R$};
\draw[color=black] (1.74,7.65) node {$L$};
\draw [fill=uuuuuu] (4.188507936507937,6.1885079365079365) circle (2.5pt);
%\draw[color=uuuuuu] (4.34,6.57) node {$S$};
\draw [fill=xdxdff] (5.386926898509582,5.430958126330731) circle (2.5pt);
%\draw[color=xdxdff] (5.54,5.87) node {$T$};
\end{scriptsize}
\end{tikzpicture}
    \caption{The $(4,4)$ grid and an extra half grid line}
    \label{fig:grid and line}
\end{figure}
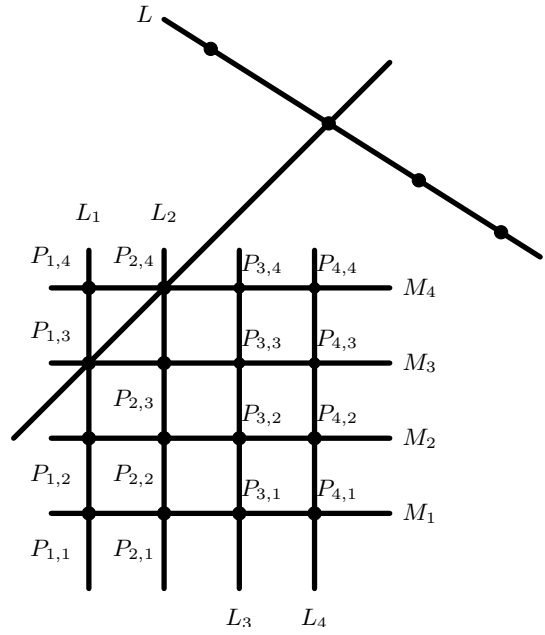

By way of an example, the permutation $\sigma^L_2$ in Figure \ref{fig:grid and line} sends $4$ to $3$.

It is convenient to consider the permutation $\sigma^L_i$ as acting on points of $Z\cap L_1$ by sending $P_{1,j}$ to $P_{1,\sigma^L_i(j)}$. This action preserves the cross-ratio of the four points by \cite[Lemma 1]{POLITUS2}, hence by Lemma \ref{lem:4 points and automorphis} it extends to a projective linear automorphism of $L_1$, which, by a slight abuse of notation, we denote with the same symbol. In the notation of \cite[Section 4]{POLITUS3} it is exactly the automorphism $f_{L_2,L_1,L}\circ f_{L_1,L_2,L_3}$.

Our assumptions impose strong conditions on the permutations $\sigma^L_i$. As already mentioned, as automorphisms of $L_1$, they preserve the cross-ratio. Additional properties are summarized  in the following lemma.
\begin{lemma}\label{lem:perm_properties}
The permutations $\sigma^L_i$ have the following properties:
\begin{enumerate}[(a)]
    \item they have no fixed points;
    \item for $i\neq j$ and for $k=1,\dots,4$ we have $\sigma^L_i(k)\neq \sigma^L_j(k)$;
%    \item each permutation $\sigma^L_i$ preserves the cross ratio $j$, in the sense that
%$$ j(P_{11},P_{12},P_{13},P_{14}) = j(P_{1\sigma^L_i(1)},P_{1\sigma^L_i(2)},P_{1\sigma^L_i(3)},P_{1\sigma^L_i(4)});$$
    \item for at least one $i\in\left\{2,3,4\right\}$ the permutation $\sigma^L_i$ is not an involution.
    \item Moreover, $\sigma_i^L$ as an automorphism of $L_1$ has exactly two fixed points; these are the intersection points of $L_1$ with the transversals.
\end{enumerate}
\end{lemma}
\begin{proof}
(a) If $P_{1j},P_{ij}$ are collinear with a point $P\in L\cap Z$, then $P=L\cap M_j$, so that $L$ intersects the quadric $Q$ in three points ($P$ and the points on the transversals), which is excluded by our assumptions. \\
(b) Assume, by contradiction, that $\sigma^L_i(k) = \sigma^L_j(k) = q$. According to (a), we know that $q \neq k$. Therefore, both lines $P_{1q},P_{ik}$ and $P_{1q},P_{jk}$ intersect with $L$. If they intersect at the same point $P$, then $M_k$ intersects $L_i$ at both $P_{iq}$ and $P_{ik}$, leading to a contradiction. Thus, the plane spanned by $P_{1q},P_{ik},P_{jk}$ intersects $L$ at two distinct points, indicating that it contains $L$. By the same reasoning, the plane also contains $M_k$, and consequently, it contains $P_{1k}$, implying that it contains $L_1$. This is impossible since $L_1$ and $L$ are skew.\\
(c) This is Lemma 6 in \cite{POLITUS2}.\\
(d) This property follows directly from the construction of $\sigma_i^L$.
\end{proof}

\subsection{Permutations and automorphism of the projective line}\label{Xing}

\begin{remark} By Lemma \ref{lem:4 points and automorphis} any permutation $\sigma$ of $P_1,P_2,P_3,P_4\in \PP^1$ which leaves the cross ratio invariant determines a unique automorphism of $\PP^1$ which restricts to $\sigma$ on the four points.
Moreover all automorphisms arising in this way have two fixed points, since any automorphism of $\PP^1$ with only one fixed point has all orbits infinite, with the exception of the fixed point.
\end{remark}

Working now with specific coordinates, we will examine which permutations from the group $S_4$ may appear as $\sigma_i^L$'s. To begin with we fix projective coordinates on $L_1$ so that
$$ P_{1,1}=[1:0],\quad P_{1,2}=[0:1],\quad P_{1,3}=[1:1], \quad P_{1,4}=[1:q]$$
with $q\neq 0,1$ and we consider the four permutations from Remark \ref{rem:CR invariant} keeping the cross  ratio of arbitrary points invariant. In Table \ref{tab:CR general} we present the associated linear maps and we determine their fixed points. Additionally we list explicitly fixed points for $q=-1$.
\renewcommand{\arraystretch}{2.9}
\begin{table}[H]
    \centering
    \begin{tabular}{c|c|c|c}
    permutation & \specialcell{linear\\ automorphism} & fixed points & $q=-1$\\
    \hline
$(2,1,4,3)$ & 
\bgroup\renewcommand*{\arraystretch}{1}$\begin{pmatrix} 
0 & 1 \\ q & 0    \end{pmatrix}$\egroup &
\specialcell{$[1:a]$ and $[1:-a]$\\ with $a^2=q$}
& $[1:\pm i]$\\

        \hline
$(3,4,1,2)$ & 
\bgroup\renewcommand*{\arraystretch}{1}$\begin{pmatrix} 
q & -1 \\ q & -q    \end{pmatrix}$\egroup &
\specialcell{$[1:q+a]$ and $[1:q-a]$\\ with $a^2=q^2-q$}
& $[1:-1\pm\sqrt{2}]$\\

    \hline
$(4,3,2,1)$ & 
\bgroup\renewcommand*{\arraystretch}{1}$\begin{pmatrix} 
1 & -1 \\ q & -1    \end{pmatrix}$\egroup &
\specialcell{$[1:1+a]$ and $[1:1-a]$\\ with $a^2=1-q$}
& $[1:1\pm\sqrt{2}]$\\
\end{tabular}
    \caption{Permutations and the associated linear maps}
    \label{tab:CR general}
\end{table}
By direct inspection we see that the fixed points of all three involutions are different. The discussion so far has the following important consequence.
\begin{corollary}\label{cor:points must be special}
Under the assumptions in the first two paragraphs of Section \ref{sec: permutations on half grids}, the points on each half grid line must be in a special position, i.e., they are either harmonic or anharmonic. 
\end{corollary}
\begin{proof}
The permutations imposed by our assumptions keep the intersection points with the transversals fixed. Since there are at least two such permutations by Lemma \ref{lem:perm_properties} (b), it is clear that the three general permutations listed in Table \ref{tab:CR general} are not enough as they all have mutually different fixed points.
\end{proof}
Every permutation in the anharmonic case (see \cite[displayed formula (3)]{POLITUS2} for an explicit list) other than one of those already considered in Table \ref{tab:CR general} has a fixed point, which is excluded by Lemma \ref{lem:perm_properties}. So we are left with the harmonic case in which we assume $q=-1$.
Among the four additional permutations there are two with fixed points, so they are excluded by Lemma \ref{lem:perm_properties} (a). The fixed points of the automorphisms generated by the remaining two permutations are presented in Table \ref{tab:CR harmonic}.
\renewcommand{\arraystretch}{2.9}
\begin{table}[H]
    \centering
    \begin{tabular}{c|c|c}
    permutation & \specialcell{linear\\ automorphism} & fixed points\\
    \hline
$(3,4,2,1)$ & 
\bgroup\renewcommand*{\arraystretch}{1}$\begin{pmatrix} 
1 & -1 \\ 1 & 1    \end{pmatrix}$\egroup &
$[1:i]$ and $[1:-i]$\\

        \hline
$(4,3,1,2)$ & 
\bgroup\renewcommand*{\arraystretch}{1}$\begin{pmatrix} 
1 & 1 \\ -1 & 1    \end{pmatrix}$\egroup &
$[1:i]$ and $[1:-i]$
\\

\end{tabular}
    \caption{The fixed points of the automorphisms generated by the 
two permutations}
    \label{tab:CR harmonic}
\end{table}
We conclude the considerations in this Section with the following result.
\begin{theorem}\label{thm:only 3 permutations}
    The only candidates for permutations $\sigma_i^L$ are
    $$(2,1,4,3),\quad (3,4,2,1)\; \mbox{ and }\; (4,3,1,2).$$
\end{theorem}
This result has important consequences which we list below.
\begin{corollary}\label{max4}
The quadric $Q$ contains at most 4 lines of the half grid in the ruling of $L_1$.
\end{corollary} 
\begin{proof}
Every line with index $i\geq 2$ induces one of the permutations in Theorem \ref{thm:only 3 permutations} on points $Z\cap L_1$. Hence $i\leq 4$.    
\end{proof}
This leads directly to the following problem.
\begin{question}\label{que:number of external lines}
What is the maximal number $n$ for which there exists a $[4,n]$-half grid?
\end{question}
We address this question in the next section for half grids with two transversals.

\section{Construction}
Now we want to apply our findings from Section \ref{sec: permutations on half grids} to provide a full classification of $[4,n]$-half grids containing a $(4,4)$-grid, under the  assumptions in the first two paragraphs of Section \ref{sec: permutations on half grids}. We will work with explicit coordinates.

To begin with, we note that any three skew lines in $\PP^3$ with coordinates $[x:y:z:w]$ can be mapped by a projective transformation to the lines: 
$$
L_1:\; \begin{cases}    y &= 0     \\  w &= 0    \end{cases},\qquad
L_2:\; \begin{cases}    x &= 0     \\  z &= 0    \end{cases},\qquad
L_3:\; \begin{cases}    y &= x     \\  w &= z    \end{cases}.
$$
These lines are contained in the quadric $Q:\; xw-yz=0$.
By further projective transformations the four harmonic points $P_{1,1},\ldots,P_{1,4}$ on $L_1$ can be normalized to $[1:0:0:0]$, $[0:0:1:0]$, $[1:0:1:0]$ and $[1:0:-1:0]$. Then the rulings on $Q$ determine the points from $Z$ on $L_2$ and $L_3$ and we obtain our initial data as:
\renewcommand{\arraystretch}{1}
\begin{equation}\label{eq:initial data} 
\begin{array}{llll}
P_{11}=[1:0:0:0],   &  P_{21}=[0:1:0:0],   & P_{31}=[1:1:0:0],         \\ 
P_{12}=[0:0:1:0],   &  P_{22}=[0:0:0:1],   & P_{32}=[0:0:1:1],         \\ 
P_{13}=[1:0:1:0],   &  P_{23}=[0:1:0:1],   & P_{33}=[1:1:1:1],         \\ 
P_{14}=[1:0:-1:0],  &  P_{24}=[0:1:0:-1],  & P_{34}=[1:1:-1:-1].        
\end{array}
\end{equation}
In the next step we need to choose the correspondence between the three permutations listed in Theorem \ref{thm:only 3 permutations} and $\sigma_2^L$, 
$\sigma_3^L$, $\sigma_4^L$. Such a choice determines the lines $L_4$ and $L$.
\begin{lemma}\label{lem:initial data}
The initial data in \eqref{eq:initial data} together with a fixed bijection 
$$\mu:\left\{ 
(2,1,4,3), (3,4,2,1), (4,3,1,2)
\right\}\to \left\{\sigma_2^L, 
\sigma_3^L, \sigma_4^L\right\}$$
determine the lines $L$ and $L_4$.
\end{lemma}
\begin{proof}
The lines $P_{1\sigma_2^L(1)}P_{21}$, $P_{1\sigma_2^L(2)}P_{22}$, $P_{1\sigma_2^L(3)}P_{23}$, $P_{1\sigma_2^L(4)}P_{24}$ determine the quadric $Q_2$ which contains $L_1,L_2,L$, while the lines $P_{1\sigma_3^L(1)}P_{31}$, $P_{1\sigma_3^L(2)}P_{32}$, $P_{1\sigma_3^L(3)}P_{33}$, $P_{1\sigma_3^L(4)}P_{34}$ determine the quadric $Q_3$ which contains $L_1,L_3,L$. The two quadrics meet in the two transversals $T_1,T_2$, the line $L_1$, and in one further line, which must be $L$.

Now $P_{41}$ must be the point of the line $M_1$, spanned by $P_{11},P_{21},P_{31}$, determined by asking that the unique quadric $Q_4$ passing through $L_1,L$, the transversals, and $P_{41}$, contains the line $P_{1\sigma_4^L}P_{41}$. In turn $P_{41}$ determines $L_4$ in the quadric $Q$.
\end{proof}
\begin{remark}
Once we know $L$ and $L_4$, the points of $Z\cap L$ and $Z\cap L_4$ are easily determined by the construction.   
\end{remark}
Following the strategy outlined in the proof of Lemma \ref{lem:initial data}, we compute equations of $L$. Of course the choice of $\sigma_2^L$ and $\sigma_3^L$ determines $\sigma_4^L$. Our results are summarized in Table \ref{tab:Ls}.
\renewcommand{\arraystretch}{1.2}
\begin{table}[H]
    \centering
    $\begin{array}{c|c|c}
      \sigma_2^L   & \sigma_3^L & \mbox{ideal of } L\\
      \hline
       (2,1,4,3)  & (3,4,2,1) & (y+z,x-w)\\
       (2,1,4,3)  & (4,3,1,2) & (y-z,x+w)\\
       (3,4,2,1)  & (2,1,4,3) & (y-z+w,x-z+2w)\\
       (3,4,2,1)  & (4,3,1,2) & (y-2z+w,x-z+w)\\
       (4,3,1,2)  & (2,1,4,3) & (y+z-w,x+z-2w)\\
       (4,3,1,2)  & (3,4,2,1) & (y+2z-w,x+z-w)
    \end{array}$
    \caption{Equations of the external line $L$}
    \label{tab:Ls}
\end{table}
The determination of $L_4$ outlined in the proof of Lemma \ref{lem:initial data} is a bit difficult to implement in  practice, so we provide an alternative approach. Having the equation of $L$, we are in the position to determine the points of $Z\cap L$. By  way of an example we do so for the data in the first row of Table \ref{tab:Ls}. For $i=1,\ldots,4$ we compute
$$R_i=P_{1\sigma_2^L}P_{2i}\cap L$$
and obtain
$$R_1=[0:1:-1:0],\;
R_2=[1:0:0:1],\;
R_3=[1:1:-1:1],
R_4=[1:-1:1:1].$$
We check directly that 
$$R_i=P_{1\sigma_3^L}P_{3i}\cap L$$
for $i=1,\ldots,4$.
This allows us to determine the points $P_{4,i}$ as the intersection points of $Q$ with the lines $P_{1,\sigma_4^L}R_i$ different from $P_{1,\sigma_4^L}$. Specifically, we obtain:
$$
P_{41}=[-1:1:0:0],\;
P_{42}=[0:0:-1:1],\;
P_{43}=[-1:1:-1:1],\;
P_{44}=[-1:1:1:-1],$$
so that the equations of $L_4$ are $x+y=0$ and $z+w=0$.

Running the same procedure for the remaining rows in Table \ref{tab:Ls}, we obtain the same line again only for the second row.

Also the lines $L_4$ match for the pairs of rows: 3 and 5, as well as 4 and 6.

Taking the $(4,4)$-grid determined this way, together with the corresponding lines $L$, we obtain in each case a set $Z$ projectively equivalent to the $F_4$ configuration.

\section{Questions}
Theorem \ref{thm:main} shows that every $[4,s]$-half grid with two transversals 
is contained in the $[4,6]$-half grid
given by the standard construction
(namely, the one given by $F_4$).
Here we show there is no $[4,s]$-half grid
with $s>6$ even if we drop the condition
on there being two transversals,
and we raise the general question of 
maximality of the half grids given
by the standard construction.

Assume $m\geq3$. The article \cite{POLITUS} constructs examples
of $[m,n]$-half grids (where $n=m+1$ if $m$ is odd and $n=m+2$ if $m$ is even) 
of $m$ points on each of $n$ lines using what it refers to as the {\it standard construction}, which we now recall. 
It starts with a $(2,2)$-grid; let $S_1,S_2,T_1,T_2$ be the grid lines, so 
$S_1$ and $S_2$ are skew, 
$T_1$ and $T_2$ are skew, and
$S_i$ and $T_j$ meet in a single point for each $i$ and $j$.

There is a linear action of ${\mathbb C}^*$ on $\PP^3$ associated to $T_1$ and $T_2$ given as follows.
The action is the identity on $T_1\cup T_2$. For each point $p\not\in T_1\cup T_2$, there is a unique line 
$L_p$ through $p$ meeting both $T_1$ and $T_2$. We can choose a coordinate system on $L_p\cong \PP^1$
such that $T_1\cap L_p$ is $[0:1]$, $T_2\cap L_p$ is $[1:0]$, and $p=[1:1]$. Then for each $u\in {\mathbb C}^*$
we set $up=[1:u]$. If we choose a coordinate system on $\PP^3$ such that $T_1: x,y=0$ and $T_2: z,w=0$,
then the action just defined has matrix 
$\begin{pmatrix}
1 & 0 & 0 & 0\\
0 & 1 & 0 & 0\\
0 & 0 & u & 0\\
0 & 0 & 0 & u\\
\end{pmatrix}$.

There is similarly a linear action of ${\mathbb C}^*$ on $\PP^3$ associated to $S_1$ and $S_2$ given 
analogously. In particular, the action is the identity on $S_1\cup S_2$. 
For each point $p\not\in S_1\cup S_2$, there is a unique line 
$L_p$ through $p$ meeting both $S_1$ and $S_2$. We can choose a coordinate system on $L_p\cong \PP^1$
such that $S_1\cap L_p$ is $[0:1]$, $S_2\cap L_p$ is $[1:0]$, and $p=[1:1]$. Then for each $u\in {\mathbb C}^*$
we set $up=[1:u]$. If we choose a coordinate system on $\PP^3$ such that $S_1$ is $x=z=0$ and $S_2$ is $y=w=0$,
then the action just defined has matrix 
$\begin{pmatrix}
1 & 0 & 0 & 0\\
0 & u & 0 & 0\\
0 & 0 & 1 & 0\\
0 & 0 & 0 & u\\
\end{pmatrix}$.

The subgroup $U_m\subset \mathbb PGL_4(\mathbb C)$ generated by the two matrices above,
where $u$ is a primitive $m$th root of 1, is isomorphic to $C_m\times C_m$, where $C_m$ is 
the multiplicative cyclic group of order $m$. The orbit of a point contained in the plane spanned by
the lines $S_i, T_j$ is contained in that plane, but the orbit of a point $p_{00}$ not contained in any of those
four planes is an $(m,m)$-grid $G$. Indeed, by appropriately scaling the variables $x,y,z,w$, the point $p_{00}$ has coordinates
$[1:1:1:1]$ and the $U_m$-orbit of $p_{00}$ consists of the points $p_{ij}=[1:u^j:u^i:u^{i+j}]$ for $0\leq i,j<m$.
Note that this set of points is an $(m,m)$-grid. To this end note that given $i$, the points 
$[1:u^j:u^i:u^{i+j}]$ for $0\leq j<m$ are collinear; denote the line containing them by $M_i$
(it is defined by $w-u^iy=u^ix-z=0$). Similarly, given $j$, the points
$[1:u^j:u^i:u^{i+j}]$ for $0\leq i<m$ are also collinear; denote the line containing them by $L_j$
(it is defined by $w-u^jz=u^jx-y=0$).
The lines $M_i$ are pair-wise skew, as are the lines $L_j$, but $M_i\cap L_j=\{p_{ij}\}$.

The question now is: what collinear sets of $m$ points can be added to $G$ to obtain a half grid of
$m$ points on $m+1$ lines. In terms of the coordinates used above, the standard construction
gives two subsets: $Y_1$, consisting of the points $[-1:0:0:u^j]$ for $0\leq j< m$, and 
$Y_2$, consisting of the points $[0:-1:u^j:0]$ for $0\leq j< m$. 
For any $m\geq3$, $G\cup Y_i$ is an $[m,m+1]$-half grid for either $i=1$ or $i=2$.
When $m$ is even, then $G\cup Y_1\cup Y_2$ is an $[m,m+2]$-half grid.

There remains the question of whether $Y_1$ and $Y_2$ are the only two subsets.
To explore this question, 
note that a necessary condition for a set $Z$ to be an $[m,r]$-half grid on $r$ lines $A_1,\ldots, A_r$, is for 
$Z\cap(A_i\cup A_j\cup A_k)$ to be a $(3,m)$-grid. 
So suppose $L$ is a line containing a set of $m$ collinear points $q_1,\ldots,q_m$ 
such that $Z=G\cup \{q_1,\ldots,q_m\}$ is an $[m,m+1]$-half grid
with half grid lines $L_0,\ldots,L_{m-1}$ and $L$. 
We will not assume that $L\cap M_i=\varnothing$ for all $i$ (although this is the case for $Y_1, Y_2$
in the standard construction). Since $L$ is not contained in the quadric containing $G$ (because $Z$ is a half grid),
and since $m\geq3$, there must be a line $M_i$ disjoint from $L$.

The lines $M_i$ and $L_j$ (for any $j$) span the plane $\Pi_{ij}$ defined by $w-u^jz-u^iy+u^{i+j}x=0$.
Pick any point $p_{ik}\in M_i$ (but not $p_{ij}$ so $k\neq j$).
Then $p_{ik}\in L_k$, so the points of $Z$ on $L, L_j$ and $L_k$ give a $(3,m)$-grid, and this grid has a transverse grid line 
$T\subset \Pi_{ij}$
through $p_{ik}$ which meets $L$ in a point $q_r$ for some $r$
(since the points $q_r$ are the points of the $(3,m)$-grid on $L$). But $L$ is skew to $L_j$, so $L$ meets $\Pi_{ij}$ in a single point, 
which thus must be the same point $q_r$ where $T$ meets $L$. This is true for each point $p_{ik}$, $k\neq j$, so
the point $L\cap \Pi_{ij}$ is a point of concurrence of $m-1$ lines where each line goes through the point $L\cap \Pi_{ij}$
and through a point of both $M_i$ and $L_j$ (but not through $p_{ij}$).

\begin{question}
Given grid lines $M_i$ and $L_j$, how many points of concurrence in the plane $\Pi_{ij}$ are there
(meaning a point $q\in\Pi_{ij}$ not on $M_i\cup L_j$ such that for each point $p_{ik}\in M_i$, $k\neq j$, the line
through $q$ and $p_{ik}$ also contains a point $p_{lj}\in L_j$)?
\end{question}

For a given $m$, this is a purely computational question.
We know there are at least two, namely $Y_1\cap \Pi_{ij}=\{[-1:0:0:u^{i+j}]\}$ and $Y_2\cap \Pi_{ij}=\{[0:-1:u^{i-j}:0]\}$,
based on the fact the standard construction gives an $[m,m+1]$-half grid.
If these are the only two for some choice of $i$ and $j$, then there are only two for each $i$ and $j$ (since
$U_m$ is a group of linear automorphisms of $\PP^3$ which acts transitively on the points $p_{ij}$).
And if there are only two, then the standard configuration with $m$ points per line is contained in no larger half grid
with $m$ points per line, even if we do not require transversals for the half grid lines.

We checked by brute force computation for $3\leq m\leq 11$ and indeed there are only two points of concurrency in these cases.
Thus the $[m,r]$-half grid given by the standard construction (with $r=m+1$ if $m$ is odd and $r=m+2$ if $m$ is even)
is contained in no $[m,s]$-half grid with $s>r$ when $3\leq m\leq 11$.

\begin{question}
Is the previous sentence true for all $m$?
\end{question}

We also pose a final question that would, if it has an affirmative answer, finish the classification of half grids with two transversals: 

\begin{question}
Let $Z$ be an $[m,n]$-half grid with two transversals where the points on each half grid line
are a single $C_m$ orbit. Must $Z$ be contained in an $[m,r]$-half grid given by the standard construction?
(So, in particular, with $r=m+1$ if $m$ is odd and $r=m+2$ if $m$ is even.)
\end{question}

\begin{paragraph}{Acknowledgements.}
This note grew out of discussions in one of the working groups active during the Workshop on Lefschetz Properties in Algebra, Geometry, Topology and Combinatorics, held in the Fields Institute, in Toronto, Canada, in the period May 14-20, 2023. It is our pleasure to thank the Institute for providing perfect working conditions and a stimulating atmosphere.
\end{paragraph}

%\bibliographystyle{abbrv}
%\bibliography{Toronto}

\newpage

\end{document}